\documentclass[12pt]{article}

\usepackage[cmtip,arrow]{xy}
\usepackage{amsmath,amssymb,enumerate,pb-diagram,pb-xy}%,pdfsync}

\usepackage[applemac]{inputenc}

\def \al{\alpha}
\def \aug{\operatorname{aug}}

\def \bs{\backslash}
\def \C{{\mathbb C}}

\def \CP{{\cal P}}

\def \coker{\operatorname{coker}}

\def \Ext{{\rm Ext}}

\def \ga{\gamma}
\def \Ga{\Gamma}
\def \H{{\mathbb H}}
\def \Hom{{\rm Hom}}

\def \mathqed{\vspace{-30pt}\\ \qed}

\def \N{{\mathbb N}}
\def \pa{{\rm par}}
\def \ph{\varphi}

\def \prf{{\bf Proof: }}
\def \PSL{{\rm PSL}}
\def \Q{{\mathbb Q}}
\def \qed{\ifmmode\eqno \square 
		\else\noproof\vskip 12pt plus 3pt minus 9pt \fi}
\def \noproof{{\unskip\nobreak\hfill\penalty50\hskip2em\hbox{}%
     \nobreak\hfill $\square$\parfillskip=0pt%
     \finalhyphendemerits=0\par}}
\def \R{{\mathbb R}}
\def \Re{\operatorname{Re}}
\def \res{{\rm res}}
\def \RH{{\rm H}}
\def \sm{\smallsetminus}
\def \SL{{\rm SL}}

\def \Z{{\mathbb Z}}
\def \({\left(}
\def \){\right)}
\def \={{\ =\ }}

\newcommand{\tto}[1]{\stackrel{#1}{\longrightarrow}} 
\newcommand{\stack}[2]{\genfrac{}{}{0pt}{}{#1} {#2}}

\newtheorem{theorem}{Theorem}[section]

\newtheorem{lemma}[theorem]{Lemma}

\newtheorem{proposition}[theorem]{Proposition}

\begin{document}

\pagestyle{myheadings} \markright{HIGHER ORDER GROUP COHOMOLOGY}

\title{Higher order group cohomology and the Eichler-Shimura map\footnote{to appear in: J. Reine Angew. Math.}}
\author{Anton Deitmar\\ \\ \\ Dedicated to Christopher Deninger}
\date{}
\maketitle

{\bf Abstract:}
Higher order group cohomology is defined and first properties are given.
Using modular symbols, an Eichler-Shimura homomorphism is constructed mapping spaces of higher order cusp forms to higher order cohomology groups.

$$ $$

\tableofcontents

\newpage
\section*{Introduction}

In the last few years, higher order modular forms have arisen in various contexts, for instance in percolation theory \cite{KZ}, Eisenstein series formed with modular symbols \cite{CDO,Gold}, or converse theorems \cite{CF,F}.
L-functions of second order forms have been studied in \cite{DKMO,IM}, Poincar\'e series attached to higher order forms have been investigated in \cite{IO}, dimensions of spaces of second order forms have been determined in \cite{DO}.

In this paper we present an approach which focuses on functorial properties of \emph{higher order cohomology groups}.
These are the derived functors of higher order invariant functors, which is a natural generalization of classical group cohomology.
It turns out that these cohomology groups can be represented as Ext-groups over the group ring.
Representing the cohomology groups as Ext-groups has the advantage that the Ext-functors produce long exact sequences out of short exact sequences plugged into either argument.
This technique is used extensively throughout the paper.

We first introduce higher order group cohomology in general.
It turns out that for finite groups or perfect groups nothing new is gained.
For Fuchsian groups we give exact sequences which allow to compute the dimensions of cohomology groups inductively.
Finally, we define the Eichler-Shimura map through modular symbols.
It is a map
$$
S_{n+2}^q(\Ga)\ \to\ \RH_{q,\pa}^1(\Ga,V_n),
$$
where $S_{n+2}^q(\Ga)$ is the space of cusp forms of even weight $n+2$ and order $q$.
The space $\RH_{q,\pa}^1(\Ga,\R)$ is the higher order group cohomology.
For $n\ge 2$ the Eichler-Shmura map is an isomorphism.
For $n=0$ it has kernel $S_{2}^{q-1}(\Ga)$ and one always has $\dim S_{n+2}^q(\Ga)=\dim \RH_{q,\pa}^1(\Ga,V_n)$.

I thank Nikolaos Diamantis and Robin Chapman for helpful comments on the contents of this paper.

\section{General groups}
Throughout, $R$ will be a commutative ring with unit.
Let $\Ga$ be a group and  $\Sigma$ a normal subgroup.
We define a sequence of functors $({\rm {\rm H}_q^0})_{q=1,2,\dots}$ from the category of $R[\Ga]$-
modules to the category of $R$-modules.
We start with $q=1$.
For an $R[\Ga]$-module $V$ let 
$$
{\rm {\rm H}_1^0}(\Ga,\Sigma,V)=V^\Ga
$$ 
be the 
usual fix-module, i.e., the set of all $v\in V$ for which $(\ga-1)v=0$.

Next suppose ${\rm {\rm H}_q^0}$ already defined, where ${\rm {\rm H}_q^0}(\Ga,\Sigma,V)$ is an $R$-
submodule of $V$. 
Let 
$
 {\rm H}_{q+1}^0(\Ga,\Sigma,V)$ 
be the set of all $v\in V$ such that $(\ga-1)v\in {\rm H}_q^0(\Ga,\Sigma,V)$ for every $\ga\in \Ga$ and $(\sigma -1)v=0$ for every $\sigma
\in\Sigma$.

If the group $\Sigma$ is clear from the context, we also write ${\rm H}_q^0(\Ga,V)$.
This will be the case later for a Fuchsian group $\Ga$, in which case we choose $\Sigma$ to be the subgroup generated by all parabolic elements of $\Ga$.

The functor ${\rm H}_q^0$ is a left-exact functor from the category of $R[\Ga]$-modules to the category of $R$-modules.
Below we will prove this fact by representing $\RH_q^0$ as a $\Hom$-functor.
We denote its right-derived functors by ${\rm H}_q^p(\Ga,\Sigma,\cdot)$.

Consider the functor $\RH^0(\Sigma,\cdot)$ from the category of $R[\Ga]$-modules to the category of $R[\Ga/\Sigma]$-modules, mapping $V$ to the $\Sigma$-fixed vectors.
Then one has $\RH_q^0(\Ga,\Sigma,\cdot)=\RH^0_q(\Ga/\Sigma,1,\RH^0(\Sigma,\cdot))$.
The functor $\RH^0(\Sigma,\cdot)$ maps injectives to injectives, so for a $\Ga$-module $V$ there is a Grothendieck spectral sequence
$$
E_2^{r,s}\=\RH_q^r(\Ga/\Sigma,1,\RH^s(\Sigma,V)),
$$
abutting to $\RH_q^{r+s}(\Ga,\Sigma,V)$.
This has interesting consequences. For instance, if $R$ is a field of characteristic zero and the virtual cohomological dimension of $\Ga/\Sigma$ is one, then there is an isomorphism
$$
\RH_q^1(\Ga,\Sigma,V)\ \cong\ \RH_q^0(\Ga/\Sigma,1,H^1(\Sigma,V)).
$$

Let $\aug:R[\Ga]\to R$ be the augmentation map given by $\aug(\sum_\ga a_\ga\ga)=\sum_\ga a_\ga$.
Let $I=\ker(\aug)$ be the augmentation ideal.
The following simple lemma will be useful.

\begin{lemma}
The ideal $I^q$ is the $R$-span of all elements of the form
$$
(\ga_1-1)\cdots (\ga_q-1),\qquad \ga_1,\dots,\ga_q\in\Ga.
$$ 
\end{lemma}

\prf
We use induction on $q$. To start with $q=1$ let $m\in I$, then
$$
m\=\sum_\ga m_\ga\ga\=\sum_\ga m_\ga\ga-\sum_\ga m_\ga\=\sum_\ga m_\ga(\ga -1).
$$
This proves the claim for $q=1$.
The induction step is clear as $I^{q+1}=I^q I$.
\qed

Let $I_\Sigma$ be the augmentation ideal of $\Sigma$.
As $\Sigma$ is normal, $R[\Ga]I_\Sigma$ is a two-sided ideal of $R[\Ga]$.
Let $J_q$ be the ideal of $R[\Ga]$ generated by the $q$-th power $I^q$ of the augmentation ideal together with $I_\Sigma$, i.e.,
$$
J_q\= I^q+R[\Ga]I_\Sigma.
$$

For short we write $A$ for the group algebra $R[\Ga]$. Let $V$ be an $A$-module.
For any ideal $J$ of $A$, we write $V^J$ for the set of all $v\in V$ with $Jv=0$.
There is a natural identification $\Hom_A(A/J,V)\cong V^J$.
One has
$$
{\rm H}_q^0(\Ga,\Sigma,V)\=V^{I^q}\cap V^{I_\Sigma}\=\Hom_A(A/J_q,V),
$$
and hence
$$
{\rm H}_q^p(\Ga,\Sigma,V)\=\Ext_A^p(A/{J_q},V).
$$

\begin{proposition}
Let $\Sigma\subset\Ga$ be an arbitrary normal subgroup.
If $\Ga$ is finite and the order $|\Ga|$ is invertible in $R$, then the 
natural injection ${\rm H}_q^0(\Ga,V)\hookrightarrow  {\rm H}_{q+1}^0(\Ga,V)$ is an 
isomorphism.
Therefore, in this case higher order group cohomology coincides with classical group cohomology.

More generally, the same conclusion holds if there is no non-trivial $|\Ga|
$-torsion in $V^\Ga$.

For arbitrary $\Ga$, the same conclusion holds if $\Ga$ coincides with its commutator subgroup $[\Ga,\Ga]$.
\end{proposition}

\prf
By induction on $q$.
Suppose first that $q=1$.
For $v\in V$ set $P(v)= \sum_{\ga\in \Ga}\ga v$.
Then $P$ is a linear map with $P^2=|\Ga|P$ and $P=|\Ga|$ on $V^\Ga$.
Further, for every $\ga\in \Ga$ one has $(\ga-1)v\in\ker P$.
Let $v\in\, \RH_2^0(\Ga,V)$.
Then $(\ga-1)v\in\ker P\cap V^\Ga=(|\Ga|-{\rm tors})\cap V^\Ga=0$.
This implies $a\in {\rm H}_1^0(\Ga,V)$.

Next assume $q>1$ and the claim proven for $q-1$. Let $v\in\, {\rm H}_q^0(\Ga,V)
$.
For $\ga\in \Ga$ one has by induction hypothesis, $(\ga-1)v\in{\rm H}_{q-1}^0
(\Ga,V)= {\rm H}_1^0(\Ga,V)$, and therefore $v\in\, \RH_2^0(\Ga,V)
= {\rm H}_1^0(\Ga,V)
$.

The last assertion is seen as follows.
For $h,g\in\Ga$, the element $ghg^{-1} {h}^{-1}-1=(gh-hg)g^{-1} {h}^{-1}=((g-1)(h-1)-(h-1)(g-1))g^{-1} {h}^{-1}$ belongs to $I^2$.
Therefore, if $\Ga= [\Ga,\Ga]$, then $I^2=I$ and hence $I^q=I$ for every $q$.
\qed

\begin{lemma}[Cocycle representation]
The module ${\rm H}_q^1(\Ga,\Sigma,V)$ is naturally isomorphic to
$$
\Hom_A(J_q,V)/\al(V),
$$
where $\al:V\to \Hom_A(J_q,V)$ is given by $\al(v)(m)=mv$.
\end{lemma}

\prf
Write $J=J_q$.
The exact sequence
$$
0\to J\to A\to A/J\to 0
$$
gives, as part of the long exact cohomology sequence of $\Ext$ in the first argument, the  exact sequence
$$
\Hom_A(A,V)\to\Hom_A(J,V)\to\Ext_A^1(A/J,V)\to
\Ext_A^1(A,V)
$$
The last term is zero, the first can be identified with $V$ and the first map is $\al$.
\qed

In the case of classical group cohomology, which is the case $q=1$, people often use the following cocycle representation for $\RH^1(\Ga,V)$.
It is the quotient $Z^1/B^1$, where $Z^1$ is the space of all maps $f:\Ga\to V$ such that $f(\ga\tau)=\ga f(\tau)+f(\ga)$ and $Z^1$ is the subspace of all $f$ of the form $f(\ga)=\ga v-v$ for some $v\in V$.
Note that for the trivial $\Ga$-module $R$ this identifies $\RH^1(\Ga,R)$ with the set $\Hom(\Ga,R)$ of group homomorphisms into the additive group of $R$.
A natural isomorphism between these two cocycle representations is given by the map $\Psi: Z^1\to\Hom_A(I,V)$,
$$
\Psi(f)(\ga-1)\=f(\ga).
$$

\begin{lemma}[Restriction]
There is a natural restriction map
$$
\res:\ {\rm H}_q^p(\Ga,\Sigma,V)\ \to\ {\rm H}^p(\Sigma,V).
$$
The restriction respects the cocycle representations.
\end{lemma}

\prf
The restriction is induced by the natural map $R[\Sigma]/I_\Sigma\hookrightarrow A/J_k$.
The fact that it maps the cocycle representation to the cocycle representation of the group cohomology, follows from the commutativity of the following diagram and the exactness of its rows.
$$
\begin{diagram}
\node{0}\arrow{e}
\node{J}\arrow{e}
	\node{R[\Ga]}\arrow{e}
		\node{R[\Ga]/J}\arrow{e}
			\node{0}\\
\node{0}\arrow{e}
\node{I_\Sigma}\arrow{n,J}\arrow{e}
	\node{R[\Sigma]}\arrow{n,J}\arrow{e}
		\node{R[\Sigma]/I_\Sigma}\arrow{n,J}\arrow{e}
			\node{0.}
\end{diagram}
$$\qed

The following lemma will be useful later.

\begin{lemma}\label{2.1}
The quotients $J_{q-1}/J_q$ and  $J_q/IJ_q$ are free $R=A/I$-modules. 
The natural map $I_\Sigma/I_\Sigma^2\to J_q/IJ_q$ is injective.
\end{lemma}

\prf
First note that $R[\Ga]$ can be identified with $R\otimes \Z[\Ga]$ and then $J_q^R=R\otimes J_q^\Z$, where we indicate the dependence of the ideal of the ring by an upper index.
Next, $J_q^\Z$ is a submodule of the free $\Z$-module $\Z[\Ga]$, hence free.
Finally, as subsets of $\Q[\Ga]$ one has $\Q J_q^\Z\cap J_{q-1}^\Z=J_q^\Z$.
This implies that $J_{q-1}^\Z/J_q^\Z$ is a free $\Z$-module and as $J_{q-1}^R/J_{q}^R\cong R\otimes J_{q-1}^\Z/J_q^\Z$, the same follows for $R$.
The proof for $J_q/IJ_q$ is analogous.

For the second assertion consider the map $I_\Sigma\to J_q/IJ_q$.
Its kernel is $I_\Sigma\cap IJ_q$, which contains $I_\Sigma^2$.
To show equality, let $f\in I_\Sigma\cap IJ_q$.
Note $IJ_q=I(I^q+ R[\Ga]I_\Sigma)=I^{q+1}+II_\Sigma$.
Accordingly, we can write $f$ as
$$
\sum_j a_j (\ga_{j,i}-1)\cdots(\ga_{j,q+1}-1)+\sum_kb_k(\ga_k-1)(\sigma_k-1),
$$
where $a_j,b_k\in R$; $\ga_{j,i}\in\Ga\sm\Sigma$; $\ga_k\in\Ga$ and $\sigma_k\in\Sigma$.
The ideal $R[\Ga]I_\Sigma$ is the kernel of the map $R[\Ga]\to R[\Ga/\Sigma]$.
It follows that the first summand is annihilated by this map and so there are $\sigma_{j,k}\in\Sigma$ and $\sigma_j\in\Sigma$ such that
$$
\sum_{j}a_j(\ga_{j,1}\sigma_{j,1}-1)\cdots(\ga_{j,q+1}\sigma_{j,q+1}-1)\sigma_j\= 0.
$$
This implies that
$$
\sum_{j}a_j(\ga_{j,1}\sigma_{j,1}-1)\cdots(\ga_{j,q+1}\sigma_{j,q+1}-1)\ \in\ R[\Ga]I_\Sigma.
$$
Now $(\ga_{j,1}\sigma_{j,1}-1)\equiv (\ga_{j,1}-1)+(\sigma_{j,1}-1)\mod II_\Sigma$, which means that $f$ lies in $II_\Sigma$.
It follows that $f$ can be written as an element of $I_\Sigma^2$ plus a sum of the form
$$
\sum_{\stack{\ga\in\Ga\sm\Sigma}{\mod\Sigma}}(\ga-1)f_\ga,
$$
where the sum runs over a fixed set of representatives and the $f_\ga$ are uniquely determined.
But $f\in R[\Sigma]$ then implies that $f_\ga=0$, so $f\in I_\Sigma^2$.
\qed

\section{Fuchsian groups}
Let $\H$ be the upper half plane in $\C$, which is acted upon via linear 
fractionals by the group $\SL_2(\R)$.
As the element $-1$ acts trivially, the action factors over $G=\PSL_2(\R)=\SL_2(\R)/\pm 1$.
Let $\Ga$ be a discrete subgroup of $G$ of finite covolume which is not cocompact.
Then $\Ga$ has finitely many equivalence classes of cusps 
$c\in\R\cup\{\infty\}$.
For each cusp $c$ fix an element $\sigma_c\in G$ such that $\sigma_c
(\infty)=c$ and $\sigma_c^{-1}\Ga_c\sigma_c=\left\{\pm\left(\begin{array}{cc}1 & 
n \\ & 1\end{array}\right): n\in\Z\right\}$, where $\Ga_c$ is the stabilizer 
subgroup of $c$ in $\Ga$.
Let $\Sigma=\Ga_\pa$ be the subgroup generated by all parabolic elements in $\Ga$.
Then $\Sigma$ is the group generated by all stabilizer groups $\Ga_c$, where $c$ varies over the cusps of $\Ga$.
Since $\ga\Ga_c\ga^{-1}=\Ga_{\ga c}$, the group $\Sigma$ is normal in $\Ga$.

Assume $\Ga$ is torsion-free.
Then there are hyberbolic generators $\ga_1,\dots,\ga_{2g}$ and parabolic generators $p_1,\dots, p_s$ such that $\Ga$ is the group generated by these with the only relation
$$
[\ga_1,\ga_2]\cdots [\ga_{2g-1},\ga_{2g}]p_1\cdots p_s\=1.
$$
The number $g\ge 0$ is called the \emph{genus} of $\Ga$.
The number $s$ is the number of inequivalent cusps of $\Ga$.
Note that the cohomological dimension of $\Ga$ is 1 if $s\ne 0$.
Otherwise it is 2 (see \cite{Schmidt}).

For $g=0$ and $q\in\N$, let $N_g(q)=0$.
For $g\ge 1$, let $N_g(q)$ be the number of all tuples $(i_1,\dots,i_q)\in\{ 1,\dots,2g\}^q$ which do not contain $(1,2)$ as a sub-tuple.
Then $N_g(1)=2g$ and $N_g(2)=(2g)^2-1$.
More generally, one has
$$
N_g(q+1)\=2gN_g(q)-N_g(q-1),
$$
which implies  that 
$$
N_g(q)\=\al^q+\al^{q-2}+\dots+\al^{-q},
$$
where $\al=g+\sqrt{g^2-1}$.
Finally, we set $N_g(0)=1$ for all $g$.

\begin{lemma}\label{dimension}
We have
$\dim J_q/J_{q+1}\=N_g(q)$.
\end{lemma}

\prf
As the case $g=0$ is trivial, we assume $g\ge 1$.
The exactness of the sequence
$$
0\to J_q/J_{q+1}\to A/J_{q+1}\to A/J_q\to 0
$$
shows that the dimension in question equals $\dim A/J_{q+1} - \dim A/J_{q}$.
In $A/J_q$ one has $p_j=1$, so we can as well assume $s=0$.
Then $J_q=I^q$.

Let $\hat A=\varprojlim_qA/I^q$ be the $I$-completion of $A$.
In $\hat A$ one has the identity
$$
\ga^{-1}\=\sum_{n=0}^\infty (1-\ga)^n,
$$
for every $\ga\in\Ga$.
So we conclude that $\hat A\cong\hat A_0$, where $A_0\subset A$ is the $R$-subalgebra generated by $\ga_1,\dots,\ga_{2g}$, further $I_0=I\cap A_0$, and $\hat A_0$ is the $I_0$-completion of $\hat A_0$.
Note that $A/I^q\cong\hat A/\hat I^q\cong \hat A_0/\hat I_0$.
It follows that $I^q/I^{q+1}$ is spanned by all elements of the form $(\ga_{i_1}-1)\cdots (\ga_{i_q}-1)$, where $1\le i_j\le 2g$ for each $j$.
The relation
$$
[\ga_1,\ga_2]\cdots [\ga_{2g-1},\ga_{2g}]\= 1
$$
is equivalent to $\ga_1\ga_2\= b\ga_2\ga_1$ where the element $b=[\ga_{2g},\ga_{2g-1}]\cdots [\ga_{4},\ga_3]$ is contained in the subgroup generated by $\ga_3,\dots\ga_{2g}$.
This means that $b$ lies in the closed subalgebra of $\hat A_0$ generated by $\ga_3,\dots\ga_{2g}$.
It follows that
$$
(\ga_1-1)(\ga_2-1)\=b(\ga_2-1)(\ga_1-1)+(b-1)(\ga_1+\ga_2-1).
$$
Modulo $I^3$, the right hand side of this equation can be written as a linear combination of elements of the form $(\ga_{i_1}-1)(\ga_{i_2}-1)$, where $(i_1,i_2)\ne (1,2)$.
As the above relation is generating all relations, there are no linear relations among the set of all elements of the form $(\ga_{i_1}-1)\cdots (\ga_{i_q}-1)$, where $(i_\nu,i_{\nu+1})\ne (1,2)$ for every $\nu=1,\dots ,q-1$.
The latter therefore form a basis of $I^q/I^{q+1}$.
The lemma follows.
\qed

Let $n\ge 0$ be an even integer.
Let $\CP_n(\R)$ be the $\R$-vector space of homogeneous polynomials 
$p(X,Y)$ of degree $n$.
So $\CP_n(\R)$ has dimension $n+1$.
There is a representation $p_n$ of $\SL_2(\R)$ on the space $\CP_n(\R)$ given 
by
$$
p_n(\ga)f\left(\begin{array}{c}X \\Y\end{array}\right)\=f\(\ga^{-1}\left(\begin{array}{c}X \\Y\end{array}\right)\).
$$
Note that the element $-1$ of $\SL_2(\R)$ acts by the scalar $(-1)^n$, so that for even $n$, one gets a representation of $G=\SL_2(\R)/\pm 1$.
From now on we set $R=\R$, and $n\ge 0$ will be an even integer.
Then the space $V=V_n=\CP_n(\R)$ is an $A=\R[\Ga]$-module, as $\Ga$ is a subgroup of $G$. 

\begin{lemma}\label{3.3}
Assume  $\Ga$ torsion-free.
Let $R=\R$.
If $n\ge 1$, then $V^{I^q}=0$, where $V=\CP_n(\R)$.
\end{lemma}

\prf
By Theorem 10.3.5 in \cite{Beard}, the group $\Ga$ contains a hyperbolic element $\ga$.
Then $\ga$ is an element of an $\R$-split torus, which implies that
 $p_n(\ga)$ is diagonalizable.
As $V$ is a highest weight module, the eigenvalue $1$ has multiplicity one and all other eigenvalues are of absolute value $\ne 1$.
Let $v\in V^{I^q}$, then for any $k_1,\dots,k_q\in\N$ we have
$$
(\ga^{k_1}-1)\cdots(\ga^{k_q}-1)v\= 0,
$$
which implies that $v$ lies in the eigenspace to the eigenvalue $1$.
As this is true for every hyperbolic element, $v$ is an element of the intersection $U$ of all $1$-eigenspaces of hyperbolic elements.
This intersection has dimension $\le 1$, so it is different from $\CP_n(\R)$.
On the other hand, $U$ is invariant under $\Ga$, and as $\Ga$ is Zariski-dense in $G$, the space $U$ is invariant under $G$.
Since the representation $p_n$ is irreducible, $U=0$.
\qed

\begin{lemma}\label{3.6}
Assume  $\Ga$ torsion-free, let $V=\CP_n(\R)$ and $q\ge 1$.
\begin{enumerate}[\rm (a)]
\item If $s\ge 1$ or $n\ge 1$, then
$\Ext_A^2(A/J_{q},V)\= 0.$
\item
If $n=s=0$, then the dimension of $\Ext_A^2(A/ J_q,\R)$ is equal to $N_g(q-1)$, and $\Ext^3(A/J_q,\R)=0$.
\end{enumerate}\end{lemma}

\prf
(a) In this case, $H^2(\Ga,V)=\Ext^2(A/I,V)$ is zero either because $\Ga$ has cohomological dimension 1 or because of Poincar\'e duality.
The exact sequence
$$
0\to J_{q}\to A\to A/J_{q}\to 0
$$
gives an isomorphism $\Ext^1_A(J_{q},V)\cong \Ext^2_A(A/J_{q},V)$.
For $q=1$ the right hand side is zero.
This implies $\Ext^1_A(I,V)=0$.
By Lemma \ref{2.1} there is an exact sequence
$$
0\to J_q\to J_{q-1}\to \R^{N}\to 0
$$
for some $N$.
As $A/I=\R$, this gives an exact sequence $\Ext_A^1(J_{q-1},V)\to \Ext_A^1(J_q,V)\to \Ext_A^2(A/I,V)^{N}=\RH^2(\Ga,V)^{N}=0$.
So the first term maps onto the second and after iteration we find a surjective map $\Ext_A^1(I,V)\to \Ext_A^1(J_q,V)$, hence both are zero.

For (b) assume $s=0$. Then $J_q=I^q$.
The cocycle representation reads $\Ext_A^1(A/I^q,\R)\cong \Hom_A(I^q,\R)\cong \Hom_\R(I^q/I^{q+1},\R)$, so the diemnsion of $\Ext^1_A(A/J_q,\R)$ equals $N_g(q)$.
Consider the exact sequence
$$
0\to J_{q-1}/J_q\to A/J_q\to A/J_{q-1}\to 0.
$$
It yields the long exact sequence
\begin{multline*}
0\to \Hom(A/J_{q-1},\R)\tto \cong\Hom(A/J_q,\R)\tto 0\Hom(J_{q-1}/J_q,\R)\\
\tto\cong\Ext^1(A/J_{q-1},\R)\tto 0\Ext^1(A/J_q,\R)\to \Ext^1(J_{q-1}/J_q,\R)\\
\to \Ext^2(A/J_{q-1},\R)\to\Ext^2(A/J_q,\R)\to\Ext^2(J_{q-1}/J_q,\R)\to 0
\end{multline*}
The first arrow is an isomorphism, so the second is zero.
The next is an isomorphism because the dimensions of the spaces match.
hence the next is zero again.
Now $\dim\Ext^1(A/J_q,\R)=N_g(q-1)$ and $\Ext^1(J_{q-1}/J_q,\R)\cong\Ext^1(\R,\R)^{N_g(q-1)}$ has dimension $2gN_g(q-1)$, since $\dim \RH^1(\Ga,\R)=2g$ by \cite{Sh} or \cite{Hida}.
Since $N_g(q)=2gN_g(q-1)-N_g(q-2)$ and $\Ext^2(J_{q-1}/J_q,\R)\cong\Ext^2(\R,\R)^{N_g(q-1)}\cong\R^{N_g(q-1)}$, the latter by Poincar\'e duality, we get an exact sequence
$$
0\to\R^{N_g(q-2)}\to\Ext^2(A/J_{q-1},\R)\to\Ext^2(A/J_q,\R)\to\R^{N_g(q-1)}\to 0.
$$
This inductively yields $\dim\Ext^2(A/J_q,\R)=N_g(q-1)$.
The last assertion is proven in a similar fashion to part (a). 
\qed

\begin{theorem}\label{3.7}
Assume  $\Ga$ torsion-free and $s>0$.
Let $R=\R$ and $V=\CP_n(\R)$.\\
\begin{enumerate}[\rm (a)]
\item If $n\ge 1$, then there is a natural exact sequence
$$
0\to {\rm H}^1_{q}(\Ga,V)\to {\rm H}_{q+1}^1(\Ga,V)\to {\rm H}^1(\Ga, V)^{N_g(q)}\to 0.
$$
\item
If  $n=0$ and $s>0$, then there is an exact sequence 
$$
0\to \R^{N_g(q)}\to \RH_{q}^1(\Ga,\R)\to {\rm H}_{q+1}^1(\Ga,\R)\to {\rm H}^1(\Ga, \R)^{N_g(q)}  \to 0.
$$
\item If $n=0=s$, then there is an exact sequence 
$$
0\to \R^{N_g(q)}\to \RH_{q}^1(\Ga,\R)\to {\rm H}_{q+1}^1(\Ga,\R)\to {\rm H}^1(\Ga, \R)^{N_g(q)}  \to \R^{N_g(q-1)}\to 0.
$$

\end{enumerate}\end{theorem}

{\bf Remark.} From \cite{Sh} or \cite{Hida} we take
$$
\dim\RH^1(\Ga,V_n)\=\begin{cases} (2g+s-2)(n+1) & n\ge 1,\\
2g+s-1 & n=0,\ s>0,\\
2g & n=0,\ s=0.\end{cases}
$$
Let $\bar N_g(q)=1+N_g(1)+\dots+ N_g(q)=\dim A/J_{q+1}$.
The theorem implies
that
$$
\dim \RH_q^1(\Ga,V_n)\=\begin{cases} \bar N_g(q-1)(2g+s-2)(n+1)& n\ge 1,\\
\bar N_g(q-1)(2g+s-2)+1 & n=0,\ s>0,\\
N_g(q) & n=0,\ s=0.
\end{cases}
$$

\prf
Consider the exact sequence
$$
0\to J_{q}/J_{q+1}\to A/J_{q+1}\to A/J_{q}\to 0,
$$
which gives 
$$
0\to \Hom(A/J_{q},V)\to\Hom(A/J_{q+1},V)\to\Hom(J_{q}/J_{q+1},V)\to
$$ $$
\to \Ext^1(A/J_{q},V)\to\Ext^1(A/J_{q+1},V)\to\Ext^1(J_{q}/J_{q+1},V)\to 0.
$$
Start with the case when $n\ge 1$.
Then the first row is zero by Lemma \ref{3.3}, and the remaining sequence reads
$$
0\to \Ext^1(A/J_{q},V)\to\Ext^1(A/J_{q+1},V)\to\Ext^1(J_{q}/J_{q+1},V)\to 0.
$$
The first assertion follows.
For the second, assume $n=0$.
Then the first map in the long exact sequence is an isomorphism and so the second is zero.
We get an exact sequence
\begin{multline*}
0\to \Hom(J_{q}/J_{q+1},\R)\to \Ext^1(A/J_{q},\R)\to\\
\to\Ext^1(A/J_{q+1},\R)\to\Ext^1(J_{q}/J_{q+1},\R)\to 0.
\end{multline*}
The last assertion is clear from this.
To prove the last item one uses the long exact sequence from the proof of Lemma \ref{3.6} (b).
\qed

\subsection{Parabolic cohomology}
Let $V$ be a $R[\Ga]$-module.
Let $C$ be the set of all cusps of $\Ga$.
For $c\in C$, let $\Ga_c$ be its stabilizer in $\Ga$.
We define the \emph{parabolic cohomology}, ${\rm H}_{q,{\rm par}}^p(\Ga,V)$ to be 
the kernel of the restriction map
$$
{\rm H}_q^p(\Ga,\Ga_\pa,V)\ \to\ \prod_{c\in C} \,{\rm H}^p(\Ga_{c},V).
$$
Any $\ga\in\Ga$ induces an isomorphism $H^p(\Ga_c,V)\to H^p(\Ga_{\ga c},V)$.
So it suffices to extend the product over a set of representatives of $\Ga\bs C$.

\begin{theorem}\label{2.5}
Assume  $\Ga$ torsion-free.
Let $R=\R$ and $V=\CP_n(\R)$.\\
\begin{enumerate}[\rm (a)]
\item If $n=0$, then there is a natural isomorphism
$$
\RH_{q,\pa}^1(\Ga,\R)\ \cong\ \Hom_\R(J_q/J_{q+1},\R).
$$
Consequently, $\dim \RH_{q,\pa}^1\= N_g(q)$.

\item If $n\ge 1$, one has an exact sequence
$$
0 \to \RH^1_{q,\pa}(\Ga,V)
\to {\rm H}_{q+1,\pa}^1(\Ga,V)\to {\rm H}_{\pa}^1(\Ga, V)^{N_g(q)}  \to 0.
$$
Consequently, $\dim \RH_{q,\pa}^1(\Ga,V)\= \bar N_g(q-1)((2g-2)(n+1)+sn)$.
\end{enumerate}\end{theorem}

\prf
(a) Let $n=0$.
The cocycle representation identifies $\RH_q^1(\Ga,\R)$ with $\Hom_A(J_q, \R)\cong\Hom_\R(J_q/IJ_q,\R)$.
Then 
$$
\RH_{q,\pa}^1(\Ga,\R)=\Hom_\R(J_q/IJ_q+I_\Sigma,\R)\cong\Hom_\R(J_q/J_{q+1},\R).
$$

(b)
Let $C$ be the set of all cusps.
We identify $\Ga\bs C$ with a set of representatives. 
Let $c$ be a cusp and $A_c=R[\Ga_c]$.
Write $I_c$ for its augmentation ideal.
Consider the exact sequence of $A$-modules,
$$
0\to J_q/J_{q+1}\to A/J_{q+1}\to A/J_q\to 0.
$$
As $A_c$-modules, these are all direct products of copies of $A_c/I_c\cong R$, hence by restriction one obtains the sequence
$$
0\to R^{N_g(q)}\to R^{\bar N_g(q)}\to R^{\bar N_g(q-1)}\to 0.
$$
We get a commutative diagram with exact rows and columns,
{\scriptsize
$$\divide\dgARROWLENGTH by2
\begin{diagram}
\node[2]{0}\arrow{s}
	\node{0}\arrow{s}
		\node{0}\arrow{s}\\
\node{0}\arrow{e}
\node{{\rm H}^1_{q,\pa}(\Ga,V)}\arrow{e}\arrow{s}
	\node{{\rm H}_{q+1,\pa}^1(\Ga,V)}\arrow{e}\arrow{s}
		\node{{\rm H}_{\pa}^1(\Ga,V)^{N_g(q)}}\arrow{s}\\
\node{0}\arrow{e}
\node{{\rm H}^1_q(\Ga,V)}\arrow{e}\arrow{s,r}{\eta_q}
	\node{{\rm H}_{q+1}^1(\Ga,V)}\arrow{e,A}\arrow{s,r}{\eta_{q+1}}
		\node{{\rm H}^1(\Ga,V)^{N_g(q)}}\arrow{s,r}\ph\\
\node{0}\arrow{e}
	\node{\prod_{c\in\Ga\bs C}\RH^1(\Ga_c,V)^{\bar N_g(q-1)}}\arrow{e}
		\node{\prod_{c\in\Ga\bs C}\RH^1(\Ga_c,V)^{\bar N_g(q)}}\arrow{e,t}{\al}
			\node{\prod_{c\in\Ga\bs C} \RH^1(\Ga_c, V)^{N_g(q)}}
\end{diagram}
$$}

Here the columns are exact by the very definition of parabolic cohomology.
By the snake lemma, the triviality of the cokernel of $\eta$ will give the desired exactness.
We prove this for $q=1$ first.
The $\Ga$-module $V$ induces a locally constant sheaf, also denoted $V$, on $\Ga\bs \H$ and $\RH^1(\Ga,V)$ equals the sheaf cohomology.
The restriction $\RH^1(\Ga,V)\to \RH^1(\Ga_c,V)$ is the restriction of sheaf cohomology to a cusp-section in $\Ga\bs \H$.
As $\Ga\bs \H$ is compact up to cusp sections, one gets an exact sequence,
$$
\RH_c^1(\Ga\bs \H,V)\to \RH^1(\Ga,V)\to \prod_{c\in\Ga\bs C}\RH^1(\Ga_s,V)\to \RH^2_c(\Ga\bs \H,V)\to 0.
$$
Here $\RH_c$ means cohomology with compact supports and
the last zero is $\RH^2(\Ga,V)=0$.
The space $\RH_c^2(\Ga\bs\H,V)$ is dual to $\RH^0(\Ga\bs\H,V)=\RH^0(\Ga,V)$ by Poincar\'e duality.
The latter space is zero, as we assume $n\ge 1$.

This argument applies to both $\eta_1$ and $\ph$, so both cokernels are zero.
By the snake lemma, we have an exact sequence,
$$
\coker(\eta_q)\to\coker(\eta_{q+1})\to\coker(\ph).
$$
For $q=1$ the left and right are zero, so the middle vanishes as well.
Next for $q=2$ the same holds true.
Inductively, we get $\coker(\eta_q)=0$ for every $q$.
Finally, the dimension formula follows from
$$
\dim \RH^1_\pa(\Ga,V_n)\= \begin{cases} (2g-2)(n+1)+sn & n\ge 1,\\
2g & n=0.
\end{cases}
$$
\mathqed

\subsection{Cusp forms}
For $k\in 2\Z$ and $f:\H\to\C$ define
$$
(f|_k\ga)(z)\=(cz+d)^{-k}f(\ga z),
$$
where $\ga=\pm\left(\begin{array}{cc}* & * \\c & d\end{array}\right)\in G$.
By linearity, we extend the definition $f|_k\sigma$ to elements $\sigma$ of 
the group ring $\R[\Ga]$.
Let $k\ge 0$ be even and let $S_k(\Ga)$ be the space of cusp forms of 
weight $k$, i.e., the complex vector space of all
\begin{itemize}
\item $f:\H\to\C$ holomorphic,
\item $f|_k(\ga-1)=0$ for every $\ga\in\Ga$,
\item for every cusp $c$ of $\Ga$, the function $(f|_k\sigma_c)(z)$ is $O
(e^{-dy})$ as $y\to +\infty$ for some $d>0$.
\end{itemize}

We sometimes also write $\ga f$ for $f|_k(\ga^{-1})$, so we can write it as a left action.

We now define cusp forms of higher order.
First let $S_k^1(\Ga)=S_k(\Ga)$, so classical cusp forms are of order $1$. Next suppose $S_k^q(\Ga)$ is already defined and let $S_k^{q+1}(\Ga)$  be the space of all functions $f$ with
\begin{itemize}
\item $f:\H\to\C$ holomorphic,
\item $f|_k(\ga -1)\in S_k^q(\Ga)$ for every $\ga\in\Ga$,
\item for every cusp $c$, $(f|_k\sigma_c)(z)=O(e^{-dy})$ as $y\to\infty$ 
for some $d>0$,
\item $f|_k(\ga-1)=0$ for every parabolic element $\ga\in\Ga$.
\end{itemize}

Note that $S_k^q(\Ga)$ is annihilated by $J_q$ with $\Sigma=\Ga_\pa$.

\begin{proposition}
For $a,b,k,l\ge 0$ we have
$$
S_k^{a+1}(\Ga) S_l^{b+1}(\Ga)\ \subset\ S_{k+l}^{a+b+1}(\Ga).
$$
So the space 
$$
S\=\bigoplus_{a,k\ge 0}S_k^{a+1}(\Ga)
$$
is a bigraded algebra.

For $f\in S_k^q(\Ga)$ one has $f|_k\ga\in S_k^q(\Ga)$.
\end{proposition}

\prf
The first assertion is done by induction on $a+b$, and the second is an easy calculation.\qed

Note \cite{Hida} that for $n\ge 0$ even,
$$
\dim S_{n+2}(\Ga)\=(2g-2)(n+1) + ns.
$$

\section{Eichler-Shimura map}
In this section we will mostly be dealing with the ring $R$ specialized to the field of real numbers $\R$.
Some of the cohomological arguments will be valid in greater generality. 
In those cases we will write $R$, whereas the use of the letter $\R$ indicates that we assume the ground ring to be $R$ here.

We define a $\CP_n(\C)$-valued differential form $\delta_n$ on $\H$ by
$$
\delta_n(z)\= (X-zY)^ndz.
$$
For a smooth function $f$ on $\H$  we set
$$
\omega(f)\= 2\pi i f(z)\delta_n(z).
$$
For any $\CP_n$-valued form $\omega$ and $\ga\in G$ let 
$$
\ga_!\omega\= p_n(\ga)\ga^{-*}\omega,
$$
where $\ga^{-*}=(\ga^{-1})^*$.
We extend $\ga\mapsto \ga_!\omega$ linearly to the group ring.

\begin{lemma}
Let $R=\R$ and $q\ge 1$.
For $f\in S_{n+2}^q(\Ga)$ one has
$$
m_!\omega(f)\=0
$$
for every $m\in J_q$.
\end{lemma}

\prf
A calculation shows that for every $\ga=\left(\begin{array}{cc}* & * \\c & 
d\end{array}\right)\in\SL_2(\R)$ one has
$$
\ga^*\delta_n(z)\=(cz+d)^{-n-2}p_n(\ga)\delta_n(z).
$$
This implies $p_n(\ga^{-1})\ga^*\omega(f)=\omega(f|\ga)$ for every $\ga\in\Ga$.
Replacing $\ga$ with $\ga^{-1}$ this means $\ga_!\omega(f)=\omega(\ga f)$.
This extends linearly to $m\in\R[\Ga]$ in place of $\ga$.
For $m\in J_q$ we have $mf=0$.
\qed

Let now $f\in S_{n+2}^q$.
For $z\in\H$ and $\ga\in\Ga$ let
$$
\ph_z(f)(\ga)\=p_n(\ga)\int_z^{\ga^{-1} z}\Re\omega(f)\ \in\ \CP_n(\R).
$$
Since $f$ is holomorphic, the integral does not depend on the path from $z$ 
to $\ga z$.
One extends the map $\ph_z(f)$ linearly to the group ring $\R[\Ga]$.
This is equivalent to the classical construction of the Eichler-Shimura isomorphism as can be found in Hida's book \cite{Hida}, where one rather uses the cocycle $\ga\mapsto \ga\ph(\ga^{-1})$, however, for a group cocycle the latter agrees with $-\ph(\ga)$.

\begin{theorem}
Suppose $\Ga$ is torsion-free and $q\ge 2$.
Let $\Sigma=\Ga_\pa$.
The map $m\mapsto\ph_z(f)(m)$ is in $\Hom_A(J_q,\CP_n(\R))$.
The induced element of the space ${\rm H}_q^1(\Ga,\Ga_\pa,\CP_n(\R))$ lies in ${\rm H}_{q,\pa}^1(\Ga,\CP_n(\R))$ and does not depend on the choice of 
$z\in\H$.

If $n\ge 1$, then the map $\ph^q$ thus defined, is an isomorphism of real vector spaces,
$$
\ph^q: S_{n+2}^q(\Ga)\ \tto\cong\ {\rm H}_{q,\rm par}^1(\Ga,\CP_n
(\R)).
$$
If $n=0$, the kernel of $\ph^q$ is $S_2^{q-1}(\Ga)\subset S_2^q(\Ga)$ so $\ph^q$ induces an injection
$$
S_2^q(\Ga)/S_2^{q-1}(\Ga)\ \hookrightarrow\ {\RH}_{q,\pa}^1(\Ga,\R).
$$
One has
$$
\dim S_2^q(\Ga)\=\dim {\RH}_{q,\pa}^1(\Ga,\R)
$$
for every $q$.
\end{theorem}

\prf
Let $\ga,\tau\in\Ga$ and compute
\begin{eqnarray*}
\ph_z(f)(\ga\tau)&=& p_n(\ga)p_n(\tau)\int_z^{\tau^{-1}\ga^{-1} z}\Re\omega(f)\\
&=& p_n(\ga)p_n(\tau)\int_z^{\tau^{-1} z}\Re\omega(f)+p_n(\ga)p_n(\tau)\int_{\tau^{-1} z}^{\tau^{-1}\ga^{-1} z}\Re\omega(f)\\
&=& p_n(\ga)\ph_z(f)(\tau)+p_n(\ga)\int_{z}^{\ga^{-1} z}\Re\tau_!\omega(f).
\end{eqnarray*}
By linearity, we can replace $\tau$ with $m\in \R[\Ga]$.
In particular, for $m\in J_q$ we get
$$
\ph_z(f)(\ga m)\=p_n(\ga)\ph_z(f)(m).
$$
This is the desired $A$-linearity.

We next show independence of $z$.
So let $z'\in\H$ and for $\ga\in \Ga$ compute
\begin{eqnarray*}
\ph_z(f)(\ga)-\ph_{z'}(f)(\ga)&=& p_n(\ga)\int_z^{\ga^{-1} z}\Re\omega(f)-p_n(\ga)\int_{z'}^{\ga^{-1} z'}\Re\omega(f)\\
&=& p_n(\ga)\int_z^{z'}\Re\omega(f)-p_n(\ga)\int_{\ga^{-1} z}^{\ga^{-1} z'}\Re\omega(f)\\
&=& p_n(\ga)\int_z^{z'}\Re\omega(f)-\int_{z}^{z'}\Re\ga_!\omega(f).
\end{eqnarray*}
Replacing $\ga$ with  $m\in J_q$ one gets,
$$
\ph_z(f)(m)-\ph_{z'}(f)(m)\= p_n(m)\int_z^{z'}\Re\omega(f).
$$
This means that the left hand side is of the form $m\mapsto mv$ for some $v$ as claimed.

We next show that $\ph$ maps to the parabolic cohomology.
Note that the exponential decay at the cusps allows to extend the 
definition of $\ph_z(f)(\ga)=p_n(\ga)\int_z^{\ga^{-1} z}\Re\omega(f)$ to the case when 
$z$ is replaced by a cusp $c$, at least if we insist that the integral 
path should be the geodesic in $\H$ from $c$ to $\ga c$.
So, if $\ga\in\Ga_c$, i.e., $\ga c=c$, then $\ph_c(f)(\ga)=0$, which 
implies that $\ph_c(f)$ zero on $\Ga_c$, so $\ph(f)$ is 
indeed parabolic.

We now show the injectivity of the map $\ph$.
Note that for $\ga\in\Ga$ we have that $(\ga -1)J_{q-1}\subset J_{q}$,
so multiplication by $(\ga -1)$ induces a map from 
$\Hom_{R[\Ga]}(R[\Ga]/J_q,V)$ to $\Hom_{R[\Ga]}(R[\Ga]/J_{q-1},V)$, which is functorial in $V$.
As cohomology is a universal $\delta$-functor, we get a map
$$
(\ga -1): \ {\rm H}_q^1(\Ga,\Sigma,V)\ \to\ {\rm H}_{q-1}^1(\Ga,\Sigma, V).
$$
We get a commutative diagram
$$
\begin{diagram}
\node{S_{n+2}^q(\Ga)}\arrow{e,t}{(\ga-1)}\arrow{s,r}{\ph^q}
	\node{S_{n+2}^{q-1}(\Ga)}\arrow{s,r}{\ph^{q-1}}\\
\node{{\rm H}_{q,\pa}^1(\Ga,V)}\arrow{e,t}{(\ga-1)}
	\node{{\rm H}_{q-1,\pa}^1(\Ga,V).}
\end{diagram}
$$
Let $f$ be in the kernel of $\ph$. By induction, $\ph^{q-1}$ is injective, so then it follows that $(\ga-1)f$ is zero for every $\ga\in\Ga$, hence $f\in S_{n+2}(\Ga)$ already.
By the case $q=1$, the map $\ph(f)\in\Hom_A(J_q,V)$ then extends to $I$, hence gives a map $I/J_q\to V$.
The image of this map $\ph(f)$ lies in $V^{J_q}$, which is zero for $n\ge 1$, so $\ph(f)=0$, so $f=0$ by the injectivity of the classical Eichler-Shimura map.
This implies the injectivity in case $n\ge 1$.

In the case $n=0$ we use induction on $q$.
For $q=1$ the claim follows from the classical Eichler-Shimura Isomorphism if we formally set $S_2^0(\Ga)=0$.
Now suppose $q\ge 2$ and the claim proven for $q-1$.
Let $f\in S_2^q(\Ga)$ be in the kernel of $\ph^q$.
Then $(\ga -1)f\in\ker\ph^{q-1}= S_2^{q-2}(\Ga)$, hence $f\in S_2^{q-1}(\Ga)$.
For the other inclusion we consider the following diagram which is commutative by the construction of the Eichler-Shimura map,
$$
\begin{diagram}
\node{S_{2}^{q-1}(\Ga)}\arrow{e,J}\arrow{s,r}{\ph^{q-1}}
	\node{S_{2}^q(\Ga)}\arrow{s,r}{\ph^q}\\
\node{\RH_{q-1,\pa}^1(\Ga,\R)}\arrow{s,r}\cong
	\node{\RH_{q,\pa}^1(\Ga,\R)}\arrow{s,r}\cong\\
\node{\Hom_A(J_{q-1}/I_\Sigma,\R)}\arrow{e,t}{\res}
	\node{\Hom_A(J_q/I_\Sigma,\R).}
\end{diagram}
$$
The restriction below is equal to zero.
Hence $S_{n+2}^{q-1}(\Ga)$ maps into the kernel of $\ph^q$.
The surjectivity of $\ph^q$ for $n\ge 1$ and the dimension formula for $n=0$ follow from our dimension formulae together with Corollary 3.13 and Theorem 4.1 of \cite{DS}.
\qed

{\small Mathematisches Institut, 
Auf der Morgenstelle 10,
72076 T\"ubingen,
Germany,
\tt deitmar@uni-tuebingen.de}

%\today


\small
\begin{thebibliography}{XXX}

\bibitem{Beard}
\bf Beardon, A.F.:
\it The geometry of discrete groups.
\rm Graduate Texts in
Mathematics, 91. Springer-Verlag, New York, 1983.

\bibitem{CDO}
\bf Chinta, G.; Diamantis, N.; O'Sullivan, C.:
\it Second order modular forms. 
\rm Acta Arith. 103 (2002), no. 3, 209--223

\bibitem{CF}
\bf Conrey, J. B.; Farmer, D. W.:
\it An extension of Hecke's converse theorem. 
\rm Internat. Math. Res. Notices 1995, no. 9, 445--463.

\bibitem{DKMO}
\bf Diamantis, N.; Knopp, M.; Mason, G.; O'Sullivan, C.:
\it $L$-functions of second-order cusp forms. 
\rm Ramanujan J. 12 (2006), no. 3, 327--347.

\bibitem{DO}
\bf Diamantis, N.; O'Sullivan, C.:
\it The dimensions of spaces of holomorphic second-order automorphic forms and their cohomology.
\rm to appear in: Transactions of the AMS.

\bibitem{DS}
\bf Diamantis, N.; Sim, D.:
\it The classification of higher-order cusp forms.
\rm to appear in: J. Reine Angew. Math.

\bibitem{F}
\bf Farmer, D.:
\it Converse theorems and second order modular forms.
\rm AMS Sectional Meeting. Salt Lake City (2002).

\bibitem{Gold}
\bf Goldfeld, D.:
\it Modular forms, elliptic curves and the $ABC$-conjecture.
\rm A panorama of number theory or the view from Baker's garden (Zürich, 1999), 128-147, Cambridge Univ. Press, Cambridge, 2002.

\bibitem{Hida}
\bf Hida, H.:
\it Elementary theory of $L$-functions and Eisenstein series. 
\rm London Mathematical Society Student Texts, 26. Cambridge University Press, Cambridge, 1993.

\bibitem{IM}
\bf Imamoglu, \"O.; Martin, Y.:
\it A converse theorem for second-order modular forms of level N.
\rm preprint.

\bibitem{IO}
\bf Imamoglu, \"O.; O'Sullivan, C.:
\it Parabolic, hyperbolic and elliptic Poincar\'e series.
\rm preprint.

\bibitem{KZ}
\bf Kleban, P.; Zagier, D.:
\it Crossing probabilities and modular forms. 
\rm J. Statist. Phys. 113 (2003), no. 3-4, 431-454
 
\bibitem{Schmidt}
\bf Schmidt, T.:
\it Rational period functions and parabolic cohomology.
\rm J. Number Theory 57 (1996), no. 1, 50--65.

\bibitem{Sh}
\bf Shimura, G.:
\it Introduction to the arithmetic theory of automorphic functions.
\rm Kano Memorial Lectures, No. 1. Publications of the Mathematical Society of Japan, No. 11. Iwanami Shoten, Publishers, Tokyo; Princeton University Press, Princeton, N.J., 1971.

\bibitem{Singerman}
\bf Singerman, D.:
\it Subgroups of Fuchsian groups and finite permutation groups. 
\rm Bull. London Math. Soc. 2, 319-323 (1970).


\end{thebibliography}
\end{document}